\documentclass[a4paper,12pt,centertags,reqno,psamsfonts]{amsart}

\usepackage[headings]{fullpage}
\usepackage{kerkis}
\usepackage{letltxmacro}
\usepackage{setspace}
\usepackage[inline]{enumitem}
\usepackage{pifont}
\usepackage{array}
\usepackage{fancybox}
\usepackage{hyperref}
\usepackage{ifthen}
\usepackage{graphicx}
\usepackage{subfigure}
\usepackage{amsmath}
\usepackage{amsthm}
\usepackage{amssymb}				
\usepackage{mathtools}
\usepackage{mathrsfs}
\usepackage{stmaryrd}
\usepackage{yhmath}
\usepackage{accents}
\usepackage{xfrac}
\usepackage[all]{xy}
\usepackage[round,comma,authoryear,sort&compress]{natbib}

\graphicspath{{Figures/}}

\setlist[itemize]{leftmargin=*,topsep=1ex,itemsep=0ex}
\setlist[enumerate]{leftmargin=*,topsep=1ex,itemsep=0ex}

\theoremstyle{plain}
\newtheorem{theorem}{Theorem}

\theoremstyle{definition}

\newtheorem{example}[theorem]{Example}

\theoremstyle{remark}

\numberwithin{theorem}{section}
\numberwithin{equation}{section}
\numberwithin{figure}{section}
\numberwithin{table}{section}

\bibpunct{(}{)}{,}{a}{}{,}

\allowdisplaybreaks[4]

\makeatletter
\let\oldr@@t\r@@t
\def\r@@t#1#2{%
\setbox0=\hbox{$\oldr@@t#1{#2\,}$}\dimen0=\ht0
\advance\dimen0-0.2\ht0
\setbox2=\hbox{\vrule height\ht0 depth -\dimen0}%
{\box0\lower0.4pt\box2}}
\LetLtxMacro{\oldsqrt}{\sqrt}
\renewcommand*{\sqrt}[2][\ ]{\oldsqrt[#1]{#2}}
\makeatother

\newcommand{\N}{\mathbb{N}}

\newcommand{\R}{\mathbb{R}}

\newcommand{\Cpq}[2]{
\ifthenelse{#1=0 \and #2=0}{\mathrm{C}}
{\ifthenelse{#2=0}{\mathrm{C}^{#1}}
{\mathrm{C}^{#1,#2}}}
}

\newcommand{\e}{\mathrm{e}}

\renewcommand{\d}{\mathrm{d}}

\newcommand{\SgAlg}[1]{\mathscr{#1}}

\newcommand{\Fltrn}[1]{\mathfrak{#1}}
\newcommand{\StpTm}{\mathfrak{S}}

\newcommand{\E}{\textsf{\upshape E}}

\renewcommand{\P}{\textsf{\upshape P}}

\newcommand{\M}{\mathscr{M}}

\newcommand{\Mloc}{\mathscr{M}_\mathrm{loc}}

\newcommand{\Hp}[1]{\mathscr{H}^{#1}}

\newcommand{\Hploc}[1]{\mathscr{H}^{#1}_\mathrm{loc}}

\newcommand{\BMO}{\mathrm{BMO}}

\begin{document}
\title{A Remark on $\Hp{1}$ Martingales}

\author{Hardy Hulley}
\author{Johannes Ruf}

\address{Hardy Hulley\\
Finance Department\\
University of Technology Sydney\\
P.O. Box 123\\
Broadway, NSW 2007\\
Australia}
\email{hardy.hulley@uts.edu.au}

\address{Johannes Ruf\\
Department of Mathematics\\
London School of Economics and Political Science\\
Houghton Street\\
London WC2A 2AE\\
United Kingdom}
\email{j.ruf@lse.ac.uk}



\date{\today}

\begin{abstract}
The space of $\Hp{1}$ martingales is interesting because of its duality with the space of $\BMO$ martingales. It is straightforward to show that every $\Hp{1}$ martingale is a uniformly integrable martingale. However, the converse is not true. That is to say, some uniformly integrable martingales are not $\Hp{1}$ martingales. This brief note provides a template for systematically constructing such processes.
\end{abstract}

\maketitle

\section{Introduction}
\label{Sec1}
Let $(\Omega,\SgAlg{F},\Fltrn{F},\P)$ be a filtered probability space, whose filtration $\Fltrn{F}=(\SgAlg{F}_t)_{t\geq 0}$ is right-continuous. By assumption, all processes are defined on $(\Omega,\SgAlg{F},\Fltrn{F},\P)$ and have c\`adl\`ag sample paths. The families of local martingales and uniformly integrable martingales are denoted by $\Mloc$ and $\M$, respectively. Given $p\geq 1$, let $\Hp{p}$ be the family of local martingales $M\in\Mloc$, for which $\E([M]_\infty^{\sfrac{p}{2}})<\infty$. Then $\Hp{p}$ forms a Banach space, when endowed with the norm $\|\cdot\|_p$, defined by
\begin{equation*}
\|M\|_p\coloneqq \Bigl(\E\bigl([M]_\infty^{\sfrac{p}{2}}\bigr)\Bigr)^{\sfrac{1}{p}},
\end{equation*}
for all $M\in\Hp{p}$. Let $\Hploc{p}$ denote the family of local martingales that are locally in $\Hp{p}$.

The space $\Hp{1}$ has several interesting features. For example,   $\Hp{1}$ contains $\Hp{2}$ as well as all local martingales with integrable variation \citep[see e.g][Theorem~IV.49]{Pro05}.  It can also be shown that $\Hp{2}$ and the family of bounded (uniformly integrable) martingales are both dense in $\Hp{1}$ \citep[see e.g][Theorem~IV.50]{Pro05}. In addition, $\Mloc\subseteq\Hploc{1}$ \citep[see e.g][Theorem~IV.51]{Pro05}. However, the most significant result is the duality between $\Hp{1}$ and the space of BMO martingales. In detail, let
\begin{equation*}
\BMO\coloneqq\bigl\{M\in\Hp{2}\,|\,\E\bigl((M_\infty-M_{\tau-})^2\bigr)\leq c^2,\;\text{for all $\tau\in\StpTm$ and some $c\in\R_+$}\bigr\},
\end{equation*}
where $\StpTm$ denotes the family of stopping times defined on $(\Omega,\SgAlg{F},\Fltrn{F},\P)$ and subject to the convention $M_{0-}\coloneqq 0$, for all $M\in\Mloc$. Then $\BMO$ becomes a Banach space, when endowed with the norm $\|\cdot\|_\BMO$, defined by
\begin{equation*}
\|M\|_\BMO\coloneqq\sup_{\tau\in\StpTm}\sqrt{\frac{\E\bigl((M_\infty-M_{\tau-})^2\bigr)}{\P(\tau<\infty)}},
\end{equation*}
for all $M\in\BMO$. It follows that $(\Hp{1})^*\simeq\BMO$ \citep[see e.g][Theorem~IV.55]{Pro05}.\footnote{This result is the probability-theoretic analogue of the classical $H^1$--BMO duality, due to \citet{FS72}.}

Now, suppose $M\in\Hp{1}$, in which case the Burkholder-Davis-Gundy inequalities \citep[see e.g][Theorem~IV.48]{Pro05} imply that $\E\bigl(\sup_{t\geq 0}|M_t|\bigr)<\infty$. An application of the dominated convergence theorem then gives 
$M\in\M$. This establishes that $\Hp{1}\subseteq\M$. However, the reverse inclusion does not hold. We demonstrate this fact below, by providing a recipe for constructing non-negative uniformly integrable martingales in $\M\setminus\Hp{1}$.
\section{Construction of Processes in $\M\setminus\Hp{1}$}
\label{Sec2}
Given a local martingale $M\in\Mloc$, observe that
\begin{equation*}
\E\biggl(\sup_{t\geq 0}|M_t|\biggr)=\int_0^\infty\P\biggl(\sup_{t\geq 0}|M_t|>u\biggr)\,\d u
\end{equation*}
and
\begin{equation*}
\sum_{n=1}^\infty\P\biggl(\sup_{t\geq 0}|M_t|>n\biggr)
\leq\int_0^\infty\P\biggl(\sup_{t\geq 0}|M_t|>u\biggr)\,\d u
\leq 1+\sum_{n=1}^\infty\P\biggl(\sup_{t\geq 0}|M_t|>n\biggr).
\end{equation*}
Since $M\in\Hp{1}$ if and only if $\E(\sup_{t\geq 0}|M_t|)<\infty$, by virtue of the Burkholder-Davis-Gundy inequalities, it follows that $M\in\Hp{1}$ if and only if
\begin{equation}
\label{eqSec1:1}
\sum_{n=1}^\infty\P\biggl(\sup_{t\geq 0}|M_t|>n\biggr)<\infty.
\end{equation}
This condition plays a key role in our construction.

Fix a non-negative local martingale $M\in\Mloc\setminus\M$ that is not a uniformly integrable martingale, and define the non-decreasing sequence $(c_n)_{n\in\N}\subset(1,\infty)$, by setting
\begin{equation}
\label{eqSec2:1}
c_n\coloneqq \ln\biggl(\e+\sum_{k=1}^n\P\biggl(\sup_{t\geq 0}M_t >k\biggr)\biggr),	
\end{equation}
for each $n\in\N$. Since $M\notin\Hp{1}$, it follows that $\lim_{n\uparrow\infty}c_n=\infty$. Next, suppose that $(\Omega,\SgAlg{F},\P)$ accommodates a discrete $\SgAlg{F}_0$-measurable random variable $Y\in\N$ that is independent of $M$, and whose distribution satisfies $\P(Y>n)=\sfrac{1}{c_n}$, for each $n\in\N$, and let
\begin{equation}
\label{eqSec2:2}
\sigma\coloneqq\inf\{t\geq 0\,|\,M_t>Y\}	
\end{equation}
denote the first time $M$ exceeds $Y$. It follows that 
\begin{equation*}
\P\biggl(\sup_{t\geq 0} M^\sigma_t>n\biggr)
\geq\P\biggl(\sup_{t\geq 0}M_t>n\biggr)\P(Y>n)
=\frac{1}{c_n}\P\biggl(\sup_{t\geq 0}M_t>n\biggr),
\end{equation*}
for each $n\in\N$. Consequently,
\begin{equation*}
\sum_{n=1}^\infty\P\biggl(\sup_{t\geq 0}M^\sigma_t>n\biggr)
\geq\lim_{m\uparrow\infty}\frac{1}{c_m}\sum_{n=1}^m\P\biggl(\sup_{t\geq 0}M_t>n\biggr)
=\lim_{m\uparrow\infty}\frac{\e^{c_m}-\e}{c_m}=\infty,
\end{equation*}
since $(c_n)_{n\in\N}$ is non-decreasing and $\lim_{n\uparrow\infty}c_n=\infty$. This implies that $M^\sigma\notin\Hp{1}$. On the other hand, the almost sure limit $M^\sigma_\infty\coloneqq M^\sigma_{\infty-}\in\R_+$ exists, since $M^\sigma$ is a non-negative local martingale, and hence also a non-negative supermartingale. Moreover,
\begin{equation*}
\begin{split}
\E(M^\sigma_\infty)&=\sum_{n=1}^\infty\E(M^\sigma_\infty\,|\,Y=n)\P(Y=n)
=\sum_{n=1}^\infty\E(M^{\tau_n}_\infty\,|\,Y=n)\P(Y=n)\\
&=\sum_{n=1}^\infty\E(M^{\tau_n}_\infty)\P(Y=n)
=\sum_{n=1}^\infty\E(M_0)\P(Y=n)
=\E(M^\sigma_0),
\end{split}
\end{equation*}
where 
\begin{equation*}
\tau_n\coloneqq\inf\{t\geq 0\,|\,M_t>n\}
\end{equation*}
for each $n \in \N$. Here, the penultimate equality follows from the fact that $M^{\tau_n}\in\M$, for each $n\in\N$. Consequently, $M^\sigma\in\M$.

\begin{example}
Consider a non-negative local martingale $M\in\Mloc$ that belongs to Class~($\mathcal{C}_0$), according to the terminology of \citet{NY06}, and suppose that $M_0=1$. In that case, $M$ is a strictly positive local martingale without any positive jumps, for which $M_\infty\coloneqq M_{\infty-}=0$. The construction above is then applicable, since $\E(M_\infty)=0<1=\E(M_0)$ implies that $M\in\Mloc\setminus\M$. Moreover, an application of Doob's maximal identity \citep[see][Lemma~2.1]{NY06} provides the following concrete representation for the non-decreasing sequence $(c_n)_{n\in\N}$, defined by \eqref{eqSec2:1}:
\begin{equation*}
c_n=\ln\biggl(\e+\sum_{k=1}^n\frac{1}{k}\biggr),	
\end{equation*}
for each $n\in\N$. It is then straightforward to see that $\lim_{n\uparrow\infty}c_n=\infty$, which is the crucial ingredient for showing that $M^\sigma\notin\Hp{1}$, where the stopping time $\sigma$ is given by \eqref{eqSec2:2}.
\end{example}

\bibliography{ProbFinBiblio}
\bibliographystyle{chicago}
\end{document}